\documentclass{amsart}
\usepackage{amsmath, amssymb, amsthm,amsfonts}
\usepackage{graphicx,color}
\usepackage{graphics}
\usepackage{comment}
\usepackage[OT2,T1]{fontenc}
\DeclareSymbolFont{cyrletters}{OT2}{wncyr}{m}{n}
\DeclareMathSymbol{\Sha}{\mathalpha}{cyrletters}{"58}

\DeclareMathOperator{\Sh}{Sh}
\DeclareMathOperator{\cond}{Cond}

\newcommand{\SHA}[1]{{\cyr X}(#1)}

\newcommand{\Q}{{\mathbb Q}}
\newcommand{\Z}{{\mathbb Z}}

\newcommand{\F}{{\mathbb F}}

\newcommand{\kro}[2]{\left( \frac{#1}{#2} \right) }

            \DeclareFontFamily{U}{wncy}{} 
            \DeclareFontShape{U}{wncy}{m}{n}{%
               <5>wncyr5%
               <6>wncyr6%
               <7>wncyr7%
               <8>wncyr8%
               <9>wncyr9%
               <10>wncyr10%
               <11>wncyr10%
               <12>wncyr6%
               <14>wncyr7%
               <17>wncyr8%
               <20>wncyr10%
               <25>wncyr10}{} 
\DeclareMathAlphabet{\cyr}{U}{wncy}{m}{n}

\begin {document}

\newtheorem{thm}{Theorem}

\newtheorem{lem}{Lemma}[section]
\newtheorem{prop}[lem]{Proposition}

\newtheorem{cor}[lem]{Corollary}

\theoremstyle{definition}

\newtheorem{ex}{Example}

\theoremstyle{remark}

\newtheorem*{ack}{Acknowledgement}

\title[On Shimura's Decomposition]{
On Shimura's Decomposition
}

\author{Soma Purkait}
\address{Mathematics Institute\\
	University of Warwick\\
	Coventry\\
	CV4 7AL \\
	United Kingdom}

\email{Soma.Purkait@warwick.ac.uk}
\date{\today}
\thanks{The author is supported by a 
an EPSRC Delivering Impact Award and an EPSRC New Directions grant}


\begin{abstract}
Let $k$ be an odd integer $\geq 3$ and $N$ a positive integer such that $4 \mid N$.
Let $\chi$ be an even Dirichlet character modulo $N$.
Shimura decomposes the space of half-integral weight cusp forms $S_{k/2}(N,\chi)$
as a direct sum
\begin{equation*}
S_{k/2}(N,\chi)=S_0(N,\chi) \oplus \bigoplus_{F} S_{k/2}(N,\chi,F),
\end{equation*}
where $F$ runs through all newforms of weight $k-1$,
level dividing $N/2$ and character $\chi^2$, 
the space $S_{k/2}(N,\chi,F)$ is the subspace 
of forms that are \lq\lq Shimura equivalent\rq\rq\ to
$F$, and
the space $S_0(N,\chi)$ is the subspace spanned by 
single-variable theta-series.
The explicit computation of this decomposition is important for practical applications
of a theorem of Waldspurger \cite{Waldspurger} relating
the critical values of $L$-functions of quadratic twists of newforms of even integral weight
to coefficients of modular forms of half-integral weight.

In this paper, we give a more precise definition of the summands $S_{k/2}(N,\chi,F)$
whilst proving that it is equivalent to Shimura's definition. We use our definition 
to give a practical algorithm for computing Shimura's decomposition, and illustrate
this with some examples.
\end{abstract}
\maketitle

\section{Introduction}

Let $F$ be a newform of even integral weight. A theorem of Waldspurger \cite{Waldspurger} 
expresses the critical value of the $\mathrm{L}$-function
of the $n$-th twist of $F$ in terms of coefficients of certain cusp forms
of half-integral weight. An example of this is the celebrated theorem of 
Tunnell \cite{Tunnell} which expresses $\mathrm{L}(E_n,1)$,
for the elliptic curve $E_n : Y^2=X^3-n^2 X$ with $n$ square-free,
in terms of coefficients of theta-series corresponding to certain positive-definite
ternary quadratic forms. As is well-known \cite[Chapter IV]{Kob},
Tunnell's Theorem gives an answer to the ancient congruent number problem,
partly conditional on the conjecture of Birch and Swinnerton-Dyer.
In explicit applications of Waldspurger's Theorem, for example Tunnell's Theorem,
it is necessary to compute, for a newform of integral weight $F$,
the space of cusp forms of half-integral weight that are \lq\lq Shimura equivalent\rq\rq\
to $F$.
 
Let $k$ be an odd integer $\geq 3$ and $N$ a positive integer such that $4 \mid N$.
Let $\chi$ be an even Dirichlet character modulo $N$.
Shimura decomposes the space of half-integral weight forms $S_{k/2}(N,\chi)$
as a direct sum
\begin{equation}\label{eqn:decomp}
S_{k/2}(N,\chi)=S_0(N,\chi) \oplus \bigoplus_{F} S_{k/2}(N,\chi,F),
\end{equation}
where $F$ runs through all newforms of weight $k-1$,
level dividing $N/2$ and character $\chi^2$, and
the space $S_{k/2}(N,\chi,F)$ is the subspace 
of forms that are \lq\lq Shimura equivalent\rq\rq\ to
$F$. 
The space $S_0(N,\chi)$ is the subspace spanned by 
single-variable theta-series.
The summands $S_{k/2}(N,\chi,F)$ occur in Waldspurger's Theorem
and their computation is necessary for explicit applications of that theorem. 

In this paper we give an algorithm for computing the 
decomposition~\eqref{eqn:decomp}. 
For this, we will give a more precise definition of
the summands $S_{k/2}(N,\chi,F)$ whilst showing
that our definition is equivalent to Shimura's definition.
In a forthcoming paper \cite{SomaIII}
we prove several Tunnell-like results with the help of our algorithm
for computing the Shimura decomposition.

The paper is organized as follows. In \cite{Shimura}, Shimura
gives a generating set for $S_0(N,\chi)$. In Section~\ref{section:S0}
we prove that this
generating set is in fact an eigenbasis. This allows us to
determine the dimension of $S_0(N,\chi)$.
For studying Shimura equivalence, we will need the theory
of Shimura lifts, which relates cusp forms of half-integral
weight to modular forms of even integral weight. We will
summarize what we need in Section~\ref{section:lifts}.
In Section~\ref{section:decomp} we give Shimura's definition
of the summands $S_{k/2}(N,\chi,F)$, our defintion
and we give a proof of the decomposition~\eqref{eqn:decomp} 
in which the summands have been redefined.
In Section~\ref{section:algo} we give our algorithm for computing the decomposition, 
prove its correctness, and remark on its practicality. 
We illustrate this practicality by giving explicit examples in Section~\ref{section:examples}.

\begin{ack} 
I would like to thank 
Professor Samir Siksek for introducing me 
to this very interesting problem and also for several helpful discussions. I would also 
like to thank the referee for many helpful remarks and suggestions.
\end{ack} 

\section{The space $S_0(N,\chi)$}\label{section:S0}

Let $N$ be a natural number such that $4 \mid N$. Let $\chi$ be an even Dirichlet character 
of modulus $N$. In this section we study the subspace $S_0(N,\chi)$ of
$S_{k/2}(N,\chi)$ which is defined as the subspace spanned by single-variable theta-series when $k=3$; 
for $k\geq 5$, we define $S_0(N,\chi)=0$.
We give Shimura's \cite{Shimura} definition of these theta-series and prove that they in fact form 
a basis of eigenforms; we therefore know the dimension of $S_0(N,\chi)$. 

Let $\psi$ be a primitive odd Dirichlet character of conductor $r_\psi$.
Let 
\[
h_\psi(z) := \sum_{m=1}^{\infty} \psi(m)mq^{m^2}.
 \] 
Shimura proves \cite[Proposition 2.2]{Shimura} that $h_\psi \in S_{3/2}(4 r_{\psi}^2,(\frac{-1}{.})\psi)$.
Consider the operator $V(t)$.
By definition, 
\[
V(t)(h_\psi)(z) = \sum_{m=1}^{\infty} \psi(m)mq^{tm^2} \in S_{3/2}\left( 4r_{\psi}^2t, \kro{-4t}{.}\psi \right).
\] 
In the literature (e.g.\ \cite[page 12]{Ono},
\cite[page 241]{Kohnen1}) these are called 
single-variable theta-functions.
Following Shimura \cite[page 478]{Shimura}, 
we define the space $S_0(N,\chi)$ to be a subspace of $S_{3/2}(N,\chi)$ spanned by
\begin{equation*}
\begin{split}
S=\{\ V(t)(h_\psi) :\ 4 r_{\psi}^2 t \mid N\ \text{and}\ &\ \text{$\psi$ is a primitive odd character of 
conductor $r_\psi$}\\
&  \text{such that}\ \chi = \kro{-4t}{.}\psi \ \}. 
\end{split}
\end{equation*}

The purpose of this section is to prove the following theorem.
\begin{thm}\label{thm:basisS0}
The set
$S$ constitutes a basis of Hecke eigenforms for $S_{0}(N,\chi)$ under Hecke operators $T_{p^2}$ for all primes $p$. 
In particular, the dimension of $S_0(N,\chi)$ is simply $\#S$.
\end{thm}
The proof of Theorem~\ref{thm:basisS0} is similar to the proof
of the corresponding result in weight $1/2$ by Serre and Stark \cite{Serre-Stark}.  We shall need a series of lemmas.
\begin{lem}\label{lem:S0}
$V(t)h_\psi$ is an eigenform for the Hecke operators $T_{p^2}$ for all primes $p$.
Indeed,
\[
T_{p^2} V(t) h_\psi=\begin{cases}
\psi(p)(1+p) V(t)h_\psi & \text{if $p \nmid 2t$} \\
\psi(p)p V(t) h_\psi  & \text{if $p \mid 2t$}.
\end{cases}
\]
\end{lem}
\begin{proof}
Let us write $V(t) h_\psi(z) = \sum_{n=1}^\infty a_n q^n$. 
Thus
\[
a_n=\begin{cases}
\psi(m) m & \text{if $n=tm^2$,} \\
0 &  \text{otherwise.}
\end{cases}
\]
Let $p$ be any prime.
Write $T_{p^2} V(t) h_\psi=\sum_{n=1}^\infty b_n q^n$. Then by \cite[Theorem 1.7]{Shimura},
\[
b_n=a_{p^2 n} + \kro{4tn}{p}\psi(p)a_n + \kro{-4t}{p}^2 \psi(p)^2 p a_{n/p^2}.  
\]
If $n/t$ is not the square of an integer, then $b_n=0$. Write $n=tm^2$. 
If $p \mid 2t$, then $b_n=a_{p^2 n}=a_{t p^2 m^2}= \psi(pm) pm$. This completes
the proof when $p \mid 2t$. 
Suppose $p \nmid 2t$. 
Then 
\begin{equation*}
\begin{split}
b_{n} &= a_{t p^2 m^2} + \kro{4t^2 m^2 }{p}\psi(p) a_{t m^2} + \kro{-4t}{p}^2 \psi(p)^2 p a_{t m^2/p^2} \\
&= a_{t p^2 m^2} + \kro{ m^2 }{p}\psi(p) a_{t m^2} +  \psi(p)^2 p a_{t m^2/p^2} \\
&=\begin{cases}
a_{t p^2 m^2} + \kro{m^2}{p}\psi(p)a_{t m^2} & \text{if $p \nmid m$}\\
a_{t p^2 m^2} + \psi^2(p)pa_{t m^2/p^2} & \text{if $p \mid m$}
\end{cases}\\
&=\psi(pm)pm +\psi(pm)m \\
&=(1+p)\psi(p)a_{t m^2}. \ 
\end{split}
\end{equation*}
Hence the lemma follows.
\end{proof}

\begin{lem}\label{cor:theoremS0}
Let $\psi_1$ and $\psi_2$ be primitive Dirichlet characters modulo 
$r_1$ and $r_2$
respectively,
and suppose $r_1\mid N$, $r_2 \mid N$. Let $\chi$ be a Dirichlet character modulo $N$ such that
$\psi_1(n) = \psi_2(n) = \chi(n)$ for all $n$ such that $(n,N)=1$.
Then $r_1 = r_2$ and $\psi_1 = \psi_2$. 
\end{lem}

\begin{proof}
The proof is immediate from \cite[Theorem 8.18]{Apostol}.
\end{proof}

\begin{proof}[Proof of Theorem~\ref{thm:basisS0}]

We will prove the theorem by showing that the elements of the set $S$ are linearly 
independent. Let $S = \{ \ V(t_i)(h_{\psi_i}) : 1\leq i \leq k \ \}$. We claim that
$t_i$ are all distinct. Suppose not. Then there exist $i$, $j$ such that
$t_i =t_j$. We know that $\chi = (\frac{-4t_i}{.})\psi_i = (\frac{-4t_j}{.})\psi_j$.
Thus, $\psi_i(n) = \psi_j(n)$ for all $(n,N)=1$. Since $\psi_i$ and $\psi_j$ are primitive,
we can apply Lemma~\ref{cor:theoremS0} to get that $\psi_i = \psi_j$ and that
$V(t_i)(h_{\psi_i}) = V(t_j)(h_{\psi_j})$. Hence the claim follows. We can assume 
that $t_1 < t_2 < \cdots < t_k$.

Now let $\alpha_{i}$ for $1\leq i \leq k$ be such that
\[
\alpha_{1}V(t_1)(h_{\psi_1}) + \alpha_{2}V(t_2)(h_{\psi_2}) + \cdots 
+\alpha_{k}V(t_k)(h_{\psi_k}) = 0.
\]
By the above equation and the $q$-expansion of $V(t_i)(h_{\psi_i})$,
it follows that 
\[
 \text{coefficient of $q^{t_1}$} = \alpha_{1}\psi_1(1) = 0.
\]
Hence $\alpha_{1} = 0$. Repeating the same argument with $t_2, t_3, \dots, t_k$, we obtain 
$\alpha_{2} = \alpha_{3} = \cdots = \alpha_{k} = 0$, 
completing the proof.
\end{proof}

\noindent{\bf Note.} 
Recall that for $k \geq 5$, we defined $S_0(N,\chi)=0$.
In the upcoming sections we will use the following
notation: 
\[
S_{k/2}^{\prime}(N,\chi):= S_0(N,\chi)^\perp;
\]
in words, the orthogonal complement to $S_0(N,\chi)$ with respect to the Petersson inner-product.
Thus, for $k \geq 5$,
\[
S_{k/2}^{\prime}(N,\chi)= 
S_{k/2}(N,\chi).
\]

\section{Shimura Lifts}\label{section:lifts}

For this section fix positive integers $k$, $N$ with $k\geq 3$ odd and $4 \mid N$.
Let $\chi$ be an even Dirichlet character of modulus $N$.
Let $N^\prime=N/2$. We will need the following theorem of Shimura.

\begin{thm} (Shimura) \label{thm:shimuracor}
Let $\lambda=(k-1)/2$. Let $f(z)=\sum_{n=1}^\infty a_n q^n \in S_{k/2}(N,\chi)$.
Let $t$ be a square-free integer and let $\psi_t$ be the
Dirichlet character modulo $tN$ defined by
\[
\psi_t(m)=\chi(m) \kro{-1}{m}^\lambda  \kro{t}{m}.
\]
Let $A_t(n)$ be the complex numbers defined by
\begin{equation}\label{eqn:shiexp}
\sum_{n=1}^\infty A_t(n) n^{-s} = \left( \sum_{i=1}^\infty \psi_t(i) i^{\lambda-1-s} \right)
\left( \sum_{j=1}^\infty a_{t j^2} j^{-s} \right).
\end{equation}
Let $\Sh_t(f)(z)=\sum_{n=1}^\infty A_t(n) q^n$.
Then
\begin{enumerate}
\item[(i)] $\Sh_t(f) \in M_{k-1}(N^\prime,\chi^2)$.
\item[(ii)] If $k \geq 5$ then $\Sh_t(f)$ is a cusp form.
\item[(iii)] If $k =3 $ and $f \in S_{3/2}^\prime(N,\chi)$
then $\Sh_t(f)$ is a cusp form. 
\item[(iv)] Suppose $f$ is an eigenform for $T_{p^2}$ for all primes $p$ and let $T_{p^2}f = \lambda_p f$. 
Then $\sum_{n=1}^\infty A_0(n) q^n \in M_{k-1}(N^\prime,\chi^2)$  where $A_0(n)$ is defined by 
\begin{equation}\label{eqn:shiexp1}
\sum_{n=1}^\infty A_0(n) n^{-s} = \prod_p (1 - \lambda_p p^{-s} + \chi(p)^2 p^{k-2-2s})^{-1}.
\end{equation}
In fact if $a_t \ne 0$ then $Sh_t(f)/a_t = \sum_{n=1}^\infty A_0(n) q^n$.
\end{enumerate} 
\end{thm}
\begin{proof}
For (i), (ii) and (iv) see \cite[Section 3, Main Theorem, Corollary]{Shimura}, for the rest see
\cite[Theorem 3.14]{Ono}. In particular, the fact that $N^\prime=N/2$ was proved by Niwa \cite[Section 3]{Niwa}.
\end{proof}
The form $\Sh_t(f)$ is called the {\em Shimura lift of $f$
corresponding to $t$}.
The following property of Shimura lifts is well known; see for example
 \cite[Chapter 3, Corollary 3.16]{Ono} or 
\cite{Kohnen}. 

\begin{prop}
Suppose $f \in S_{k/2}(N,\chi)$.  Let $t$ be a square-free positive integer. If $p \nmid tN$
is a prime then
\[
\Sh_t( T_{p^2} f)=T_{p} \Sh_t(f). 
\]
\end{prop}
Here $T_{p^2}$ is the Hecke operator on $S_{k/2}(N,\chi)$ and $T_p$ is the 
Hecke operator on $M_{k-1}(N^\prime,\chi^2)$.

In \cite{thesis}, we prove the following strengthening of this result.
\begin{prop}\label{prop:commute}
Suppose $f \in S_{k/2}(N,\chi)$ and $t$ a square-free
positive integer.  If $p$
is a prime then
\[
\Sh_t( T_{p^2} f)=T_{p} \Sh_t(f). 
\]
\end{prop}

In a forthcoming paper \cite{SomaII} we will be using this stronger result to give the generators 
of Hecke algebra as a $\Z[\zeta]$-module, where $\zeta$ is a primitive
$\varphi(N)$-th root of unity; here $\varphi$ stands for the Euler's totient function. 

\section{Shimura's Decomposition}\label{section:decomp}

In this section we state and refine a theorem of Shimura
that conveniently decomposes the space of cusp forms 
of half-integral weight.

As before let $k$, $N$ be positive integers with $k \ge 3$ odd and $4 \mid N$. Let 
$\chi$ be an even Dirichlet character of modulus $N$.
Let $N^\prime=N/2$. For $M \mid N^\prime$ such that $\cond(\chi^2) \mid M$
and 
a newform $F \in S_{k-1}^{\mathrm{new}}(M,\chi^2)$
define
\[
S_{k/2}(N,\chi,F)=
\{
f \in S_{k/2}^{\prime}(N,\chi) : 
\text{$T_{p^2}(f)=\lambda_{p}^{F} f$ for almost all $p \nmid N$}\}; 
\]
here $T_{p} (F) =\lambda_{p}^{F} F$.
\begin{thm}\label{thm:Shimura} (Shimura \cite{Shimura2}) We have
$S_{k/2}^{\prime}(N,\chi)=\bigoplus_F S_{k/2}(N,\chi,F)$ where $F$
runs through all newforms $F \in S_{k-1}^{\mathrm{new}}(M,\chi^2)$
with $M \mid N^\prime$ and $\cond(\chi^2) \mid M$.
\end{thm}

For us this theorem is not suitable for computation since for
any particular prime $p \nmid N$, we do not know if it is included
or excluded in the \lq almost all\rq\ condition. In fact we shall prove this
theorem with a more precise definition for the spaces $S_{k/2}(N,\chi,F)$.

From now on and for the rest of the paper we take the following as
the definition of the space $S_{k/2}(N,\chi,F)$:
\[
S_{k/2}(N,\chi,F)=
\{
f \in S_{k/2}^{\prime}(N,\chi) : 
\text{$T_{p^2}(f)=\lambda_{p}^{F} f$ for  all $p \nmid N$}\}. 
\]

We say that $f \in S_{k/2}^{\prime}(N,\chi)$ is {\em Shimura equivalent} to $F$ if $f$ belongs to the space 
$S_{k/2}(N,\chi,F)$.

\begin{thm}\label{thm:decomp}
Shimura's decomposition in Theorem~\ref{thm:Shimura}
holds with this new definition.
\end{thm}

\begin{proof}
It is well-known that the operators $ \overline{\chi(p)} T_{p^2}$ 
on $S_{k/2}(N,\chi)$ 
with $p \nmid N$ commute and are Hermitian; see for example \cite{Kob}. 
They also preserve the subspace $S_{k/2}^{\prime}(N,\chi)$.
Therefore, there exists
an eigenbasis $f_1,f_2,\dots,f_n$ for $S_{k/2}^{\prime}(N,\chi)$
with respect to the operators $T_{p^2}$ for $p \nmid N$.
Let $f$ be one of the $f_i$. Let $H=\Sh_t(f)$, the Shimura lift of $f$ corresponding to $t$ 
(Theorem~\ref{thm:shimuracor}) with any square-free $t$. 
We know that $H \in S_{k-1}(N^\prime,\chi^2)$.
Moreover, for all $p \nmid N$ we know that $H$ is an eigenfunction
for $T_p$ and it has the same eigenvalue as $f$ under $T_{p^2}$;
see Proposition~\ref{prop:commute}. By the theory of newforms of integral weight modular forms 
\cite[Corollary 4.6.20]{Miyake} we know that there exists
uniquely a divisor $M$ of $N^\prime$ with $\cond(\chi^2) \mid M$ and a newform
$F \in S_{k-1}^{\mathrm{new}}(M,\chi^2)$ such that $F$
has the same $T_p$-eigenvalues as $H$ for all primes $p \nmid N^\prime$.
Thus $f \in S_{k/2}(N,\chi,F)$. Hence $S_{k/2}^{\prime}(N,\chi)$ is a sum of the subspaces $S_{k/2}(N,\chi,F)$ as 
$F$ runs through newforms $F \in S_{k-1}^{\mathrm{new}}(M,\chi^2)$ with $M \mid N^\prime$ and $\cond(\chi^2) \mid M$.
We now show that this sum is actually a direct sum. For this, we just need to show that if $h_1,h_2,\dots,h_r$ are
all the elements of the above eigenbasis that belong to $S_{k/2}(N,\chi,F_0)$
where $F_0$ is a fixed newform in $S_{k-1}^{\mathrm{new}}(M_0,\chi^2)$ with $M_0 \mid N^\prime$ 
and $\cond(\chi^2) \mid M_0$, then they actually form a basis for the space $S_{k/2}(N,\chi,F_0)$.
We can reorder our basis elements such that $f_i=h_i$ for $1 \le i \le r$.
Let $h \in S_{k/2}(N,\chi,F_0)$ and suppose $h=\alpha_1f_1+\alpha_2f_2+\cdots+\alpha_nf_n$.
We show that $\alpha_i =0$ for $r+1 \le i \le n$. We will show that $\alpha_{r+1} =0$ 
and the same argument follows for the others. We know that $f_{r+1} \in S_{k/2}(N,\chi,F)$
for some suitable newform $F$ and $F_0 \ne F$. This implies there exists a prime $p \nmid N$ such that
$\lambda_p^{0} \ne \lambda_p$ where $\lambda_p^{0}$ and $\lambda_p$ are corresponding $T_p$-eigenvalues
of $F_0$ and $F$. Applying $T_{p^2}$ to $h$ we get $\alpha_{r+1} =0$.
The theorem follows.
\end{proof}

In fact, as a corollary to the proof of Theorem~\ref{thm:decomp}
we can deduce the following precise relationship between
the Shimura lift $H$ and the newform $F$.
\begin{cor}\label{cor:lincomb}
Let $F$ be a newform belonging to $S_{k-1}^{\mathrm{new}}(M,\chi^2)$
where $M \mid N^\prime$ and $\cond(\chi^2) \mid M$.
Let $f \in S_{k/2}(N,\chi,F)$ and let $H=\Sh_t(f)$ for any square-free $t$. Then we can
write $H$ as a linear combination
\[
H=\sum_{d \mid (N^\prime/M)} \alpha_d V_d(F).
\]
\end{cor}

The following is an easy lemma that shows equivalence of the two definitions thereby 
leading to an alternate proof of Theorem~\ref{thm:Shimura}. 
\begin{lem}\label{lem:agree}
Our definition of $S_{k/2}(N,\chi,F)$ agrees with Shimura's
definition. In other words, 
if we write
\[
S^{\mathrm Sh}_{k/2}(N,\chi,F)=
\{
f \in S_{k/2}^{\prime}(N,\chi) : 
\text{$T_{p^2}(f)=\lambda_{p}^{F} f$ for almost all $p \nmid N$}\};
\]
then 
$S^{\mathrm Sh}_{k/2}(N,\chi,F)= S_{k/2}(N,\chi,F)$.
\end{lem}
\begin{proof}
Clearly, the right-hand side is contained in the left-hand side.  
Suppose $f$ is in left-hand side. 
We use the decomposition Theorem~\ref{thm:decomp} with our definition 
of summands. Let $G$ run through the newforms of 
levels dividing $N/2$. Then we can write
$f=\sum f_G$ where $f_G \in S_{k/2}(N,\chi,G)$.
Here $F$ is one of the $G$s.
We know that for almost all primes $p$, 
\[
T_{p^2}f=\lambda^F_p f= \sum \lambda^F_p f_G
\]
where $T_p F = \lambda_p^F F$.
But,
\[
T_{p^2}f=\sum T_{p^2} (f_G)=\sum \lambda^G_p f_G
\]
where $T_p G=\lambda_p^G G$. Thus
\[
\sum (\lambda_p^F-\lambda_p^G) f_G=0.
\]
By the fact that the summands belong to a direct sum,
we see that each summand must individually be zero.
If $f_G \ne 0$ then $\lambda_p^F=\lambda_p^G$
for almost all $p$ which forces $G=F$ by \cite[Theorem 4.6.19]{Miyake}.
Thus $f=f_F \in S_{k/2}(N,\chi,F)$ as required.
\end{proof}

\section{Algorithm for computing Shimura's decomposition}\label{section:algo}

The following theorem gives our algorithm for computing the Shimura decomposition.
\begin{thm} \label{thm:algo}
Let $F_1,\dots,F_m$ be  
the newforms of weight $k-1$, character $\chi^2$
and level dividing $N^\prime$.
For prime $p$, 
and $F$ one of these newforms,
write $T_p (F)= \lambda_{p}^{F} F$.
Let $p_1,\dots,p_n \nmid N$ be primes such that the $m$ vectors of eigenvalues
$(\lambda_{p_1}^{F},\dots,\lambda_{p_n}^{F})$, with $F=F_1,\dots,F_m$, are pairwise distinct.
If $f \in S_{k/2}^\prime (N,\chi)$ is an eigenform for $T_{p_i^2}$
for $i=1,\dots,n$ then 
$f$ belongs to one of the summands $S_{k/2}(N,\chi,F)$.
\end{thm}
\begin{proof}
Suppose $f \in S_{k/2}^\prime(N,\chi)$ is an eigenform for $T_{p_i^2}$
for $i=1,\dots,n$. 
Write $T_{p_i^2} f = \mu_i f$.
By Shimura's decomposition, we can write
\[
f=\sum_{F} f_F
\]
for some unique $f_F \in S_{k/2}(N,\chi,F)$; here $F$ varies over $F_i$, $1\le i\le m$.  Thus
\[
\sum_{F} \lambda_{p_i}^{F} f_F=T_{p_i^2} f=\mu_i \sum_F f_F.
\]
As the decomposition is a direct sum, we find that
\[
(\lambda_{p_i}^{F}-\mu_i) f_F=0, \qquad i=1,\dots,n.
\]
We will show that at most one $f_F$ is non-zero. This will
force $f$ to be in one of the components $S_{k/2}(N,\chi,F)$
which is what we want to prove. Suppose therefore that
$f_{F_1}\ne 0$ and $f_{F_2}\ne 0$. Then
\[
\lambda_{p_i}^{F_1}=\mu_i=\lambda_{p_i}^{F_2}, \qquad i=1,2,\dots,n.
\]
This contradicts the assumption that the vectors of eigenvalues are distinct,
and completes the proof.
\end{proof}

We can reframe Theorem~\ref{thm:algo} as follows.

\begin{cor}\label{cor:prev}
Let $F$ be a newform of weight $k-1$, level $M$ dividing $N^\prime$, and character $\chi^2$.
Let $p_1,\dots,p_n$ be primes not dividing $N$ satisfying the following: for every newform 
$F^\prime \ne F$ of weight $k-1$, level dividing $N^\prime$ and character $\chi^2$, there
is some $p_i$ such that $\lambda_{p_i}^{F^\prime} \ne \lambda_{p_i}^{F}$, where
$T_{p_i}(F)=\lambda_{p_i}^{F} \cdot F$. Then
\[
S_{k/2}(N,\chi,F)=
\left\{
f \in S_{k/2}^\prime(N,\chi) \; : \;
T_{p_i^2} (f) = \lambda_{p_i}^{F} f  \quad \text{for $i=1,\dots,n$} 
\right\}.
\]
\end{cor}

Recall that $S_{k/2}^\prime(N,\chi) = S_{k/2}(N,\chi)$ except possibly when $k=3$. 
We have the following refinement of the above corollary which takes care of the case 
when $S_{k/2}^\prime(N,\chi) \subsetneq S_{k/2}(N,\chi)$, that is, $S_0(N,\chi) \ne 0$.

\begin{cor} \label{cor:algo}
Assuming the notation in the above corollary, the following stronger statement holds :
\[
S_{k/2}(N,\chi,F)=
\left\{
f \in S_{k/2}(N,\chi) \; : \;
T_{p_i^2} (f) = \lambda_{p_i}^{F} f  \quad \text{for $i=1,\dots,n$} 
\right\}.
\]
\end{cor}

\begin{proof}
Let $f_1,\dots f_r$ be the basis of eigenforms for $S_0(N,\chi)$ as stated in Theorem~\ref{thm:basisS0}. Recall that 
$f_i = V(t_i)h_{\psi_i}$ where $\psi_i$ is primitive odd character of conductor $r_{\psi_i}$ such that 
$4r_{\psi_i}^2t_i \mid N$ and $\chi = \kro{-4t_i}{.}\psi_i$. Let $q = p_i$ for some fixed $i$. We claim that 
$T_{q^2}(f_i) \ne \lambda_{q}^{F} f_i$ for any $1 \le i \le r$. Since $F$ is a newform of weight $2$ we know 
by Deligne's work on Weil conjectures that $\lvert  \lambda_{q}^{F} \rvert \le 2\sqrt{q}$. By Lemma~\ref{lem:S0}, 
$T_{q^2}(f_i) = \psi_i(q)(1 + q) f_i$ as $q \nmid N$. Clearly 
$\lvert \psi_i(q)(1 + q) \rvert = \lvert 1+q \rvert > 2\sqrt{q}$. Hence the claim follows. 

Let $g \in S_{k/2}(N,\chi)$ such that $T_{p_i^2} (g) = \lambda_{p_i}^{F} g$ for $1 \le i \le n$. We can write $g$ uniquely as 
$g = g_1 + g_2$ where $g_1 \in S_0(N,\chi)$ and $g_2 \in  S_{k/2}^\prime(N,\chi)$. Since the Hecke operators $T_{p_i^2}$ preserve the 
subspaces $S_0(N,\chi)$ and $S_{k/2}^\prime(N,\chi)$ we obtain $T_{p_i^2} (g_j) = \lambda_{p_i}^{F} g_j$ for all $1 \le i \le n$ and 
$j=1,2$. Thus by Corollary~\ref{cor:prev}, $g_2 \in S_{k/2}(N,\chi,F)$. We show that $g_1=0$. Let $g_1 = \sum_{i=1}^{r}a_if_i$. 
In particular for the prime $q$ we must have $a_iT_{q^2}(f_i) = a_i\lambda_{q}^{F} f_i$. The above claim implies that 
$a_i =0$ for all $1 \le i \le r$. Hence we are done.
\end{proof}

\noindent {\bf A remark on the practicality of our algorithm.} The working of the algorithm is based on multiplicity-one 
theorem~\cite[Theorem 4.6.19]{Miyake} of newforms of integral weight and Sturm's bound \cite{Sturm}. Indeed multiplicity-one 
guarantees existence of the primes $p_i$ in the algorithm. However thanks to Sturm's result we need to look for such primes 
only up to the Sturm's bound. In practice, the set of these primes is usually very small. For example, for the decomposition 
of the space $S_{3/2}^\prime(1984,\chi_0)$ in Example~\ref{ex:1984} below, we only need to work with primes in the set
$\{3, 5, 7, 13, 19\}$. In fact running our algorithm for level $4N$ with $4 \le 4N \le 3000$, we observe that 
$31$ is the largest prime we need to work with in order to decompose the space $S_{3/2}(4N, \chi)$ into Shimura equivalent spaces; 
here $\chi$ is a quadratic Dirichlet character modulo $4N$. Futhermore we have to go as far as the prime $31$ in only five instances. 
One could heuristically argue that given two distinct rational newforms of weight $2$ the probability that $P$ is the first 
prime where the Fourier coefficients differ is roughly $(1-\frac{1}{4\sqrt{P}}) \cdot \prod \frac{1}{4\sqrt{p}}$ where the product 
runs over primes $p$ less than $P$ with $p$ not dividing the levels of the newforms. One can develop such a heuristic argument to 
explain why the set of primes that we need for our algorithm is rather small, although it does not seem possible to
supply a rigorous proof that the set of primes is small.

\section{Examples}\label{section:examples}

\begin{ex}\label{ex:32}
 Let $\chi_0$ be the trivial Dirichlet character modulo $32$. 
 One can see using dimension formula \cite{C-O} that the space $S_{7/2}(32,\chi_0)$ has dimension $6$.
 We use Theorem \ref{thm:decomp} to obtain the following decomposition
\begin{equation*}
S_{7/2}(32,\chi_0) = S_{7/2}(32,\chi_0,G_4) \oplus S_{7/2}(32,\chi_0,G_8) \oplus S_{7/2}(32,\chi_0,G_{16})
\oplus S_{7/2}(32,\chi_0,G_{16}^{\prime})
\end{equation*}
where $G_4$, $G_8$ and $G_{16}$, $G_{16}^{\prime}$ are the newforms of weight $6$, trivial character and levels 
$4$, $8$ and $16$ respectively and are given by following $q$-expansions: 
\begin{equation*}
\begin{split}
G_4 &= q - 12q^3 + 54q^5 - 88q^7 - 99q^9 + 540q^{11} + O(q^{12})\\
G_8 &= q + 20q^3 - 74q^5 - 24q^7 + 157q^9 + 124q^{11} + O(q^{12}) \\
G_{16} &= q + 12q^3 + 54q^5 + 88q^7 - 99q^9 - 540q^{11} + O(q^{12})\\
G_{16}^{\prime} &= q - 20q^3 - 74q^5 + 24q^7 + 157q^9 - 124q^{11} + O(q^{12}).
\end{split}
\end{equation*}
We use our algorithm (Theorem \ref{thm:algo}) to obtain 
\[
S_{7/2}(32,\chi_0,G_4) = \langle g_1, g_2, g_3, g_4 \rangle, \quad S_{7/2}(32,\chi_0,G_8) = \langle g_5, g_6\rangle,\]
\[
S_{7/2}(32,\chi_0,G_{16}) = \{0\} = S_{7/2}(32,\chi_0,G_{16}^{\prime}).
\]
where $g_i$ have the following $q$-expansions:
\begin{equation*}
\begin{split}
g_1 &=  q - 3q^9 - 8q^{17} + 29q^{25} + O(q^{30})\\ 
g_2 &= q^2 - 6q^6 + 10q^{10} + 4q^{14} - 21q^{18} + 10q^{22} - 18q^{26} + O(q^{30})\\
g_3 &= q^4 - 2q^8 - 4q^{20} + 12q^{24} + O(q^{30})\\
g_4 &= q^5 - 7q^{13} + 18q^{21} - 21q^{29} + O(q^{30})\\
g_5 &= q^2 + 2q^6 - 6q^{10} - 12q^{14} + 11q^{18} + 18q^{22} - 2q^{26} + O(q^{30})\\
g_6 &= q^3 - 5q^{11} + 3q^{19} + 20q^{27} + O(q^{30})
\end{split}
\end{equation*}
\end{ex}

\begin{ex}\label{ex:72}
Let $\chi_3$ be Dirichlet character modulo $72$ given by $\chi_3(\cdot)= \kro{3}{\cdot}$. 
As before using dimension formula \cite{C-O} we can see that the space $S_{5/2}(72,\chi_3)$ has dimension $12$.
We use Theorem \ref{thm:decomp} to obtain the following decomposition
\begin{equation*}
\begin{split}
S_{5/2}(72,\chi_3) = & S_{5/2}(72,\chi_3,H_6) \oplus S_{5/2}(72,\chi_3,H_9) \oplus S_{5/2}(72,\chi_3,H_{12}) \\
& \oplus S_{5/2}(72,\chi_3,H_{18}) \oplus S_{5/2}(72,\chi_3,H_{36})
\end{split}
\end{equation*}
where $H_6$, $H_9$, $H_{12}$, $H_{18}$, $H_{36}$ are the unique newforms of weight $4$, trivial character and levels 
$6$, $9$, $12$, $18$ and $36$ respectively and are given by following $q$-expansions:
\begin{equation*}
\begin{split}
H_6 & = q - 2q^2 - 3q^3 + 4q^4 + 6q^5 + 6q^6 - 16q^7 - 8q^8 + 9q^9 - 12q^{10} + 12q^{11} + O(q^{12})\\
H_9 & = q - 8q^4 + 20q^7 + O(q^{12})\\
H_{12} & = q + 3q^3 - 18q^5 + 8q^7 + 9q^9 + 36q^{11} + O(q^{12}), \\
H_{18} & = q + 2q^2 + 4q^4 - 6q^5 - 16q^7 + 8q^8 - 12q^{10} - 12q^{11} + O(q^{12}),\\
H_{36} & = q + 18q^5 + 8q^7 - 36q^{11} + O(q^{12}).
\end{split}
\end{equation*}
We use our algorithm (Theorem \ref{thm:algo}) to obtain 
\[
S_{5/2}(72,\chi_3,H_6) = \langle h_1, h_2, h_3, h_4, h_5, h_6 \rangle, \quad S_{5/2}(72,\chi_3,H_9) = \{0\}, \]
\[
S_{5/2}(72,\chi_3,H_{12}) = \langle h_7, h_8, h_9 \rangle, \;
S_{5/2}(72,\chi_3,H_{18}) = \langle h_{10}, h_{11}\rangle, \; S_{5/2}(72,\chi_3,H_{36})= \langle h_{12} \rangle,
\]
where $h_i$ have the following $q$-expansions: 
\begin{equation*}
\begin{split}
h_1 &= q + 4q^{10} - 8q^{13} - 8q^{22} + 11q^{25} + O(q^{30})\\ 
h_2 &= q^2 - q^5 - 2q^{14} + q^{17} + 6q^{26} - 3q^{29} + O(q^{30})\\
h_3 &= q^3 - 2q^{12} - 3q^{27} + O(q^{30})\\
h_4 &= q^4 - 2q^{16} - 2q^{19} + O(q^{30})\\
h_5 &= q^8 - q^{11} - q^{20} + O(q^{30})\\
h_6 &= q^9 - 2q^{18} - 2q^{21} + O(q^{30})\\
h_7 &= q - 2q^{10} + 4q^{13} - 8q^{22} - 13q^{25} + O(q^{30})\\
h_8 &= q^2 - 4q^5 + 10q^{14} - 2q^{17} - 18q^{26} + 12q^{29} + O(q^{30})\\
h_9 &= q^6 - q^9 - q^{18} + O(q^{30})\\
h_{10} &= q - 8q^{10} + 4q^{13} + 16q^{22} - q^{25} + O(q^{30})\\
h_{11} &= q^4 - 2q^7 + 2q^{16} - 4q^{28} + O(q^{30})\\
h_{12} &= q^2 + 2q^5 - 2q^{14} - 8q^{17} - 6q^{26} + 6q^{29} + O(q^{30})
\end{split}
\end{equation*}
\end{ex}

\noindent {\bf Remark.} Given a newform $F$ of integral weight $k-1$ and level $N$, it is natural to ask the minimum
level at which one can find the Shimura equivalent forms of weight $k/2$ corresponding to $F$. This has been answered 
by Mao~\cite[Theorem 1.1]{Mao} when the level of $F$ is odd. We are interested in looking at the cases when 
$N$ is even and in particular when either $4\parallel N$ or $8\parallel N$. We note that in these particular cases 
Waldspurger's Theorem~\cite[Th\'{e}or\`{e}me 1]{Waldspurger} is not applicable. We apply our algorithm to several such examples of 
newforms. We observe that if $F$ is a newform in $S_{k-1}^{\mathrm{new}}(4N)$ with $N$ odd and square-free
then the smallest level at which there is a non-zero Shimura equivalent form is $8N$ and the space $S_{k/2}(8N, \chi_0, F)$ 
is one-dimensional; here $\chi_0$ is the trivial character. Further if $F$ is a newform in 
$S_{k-1}^{\mathrm{new}}(8N)$ with $N$ odd and square-free then the smallest level at which there is a non-zero Shimura equivalent 
form is $32N$ and the space $S_{k/2}(32N, \chi_0, F)$ is now two-dimensional. In particular, in Example~\ref{ex:32}, we find that 
$S_{7/2}(16,\chi_0,G_8) = {0}$ while  $S_{7/2}(32,\chi_0,G_8)$ is two-dimensional. Also in Example~\ref{ex:72} we find 
that for $H_{12}$ and $H_{36}$ the respective Shimura equivalent spaces $S_{5/2}(24,\chi_3,H_{12})$ and $S_{5/2}(72,\chi_3,H_{36})$ 
are each one-dimensional.

\begin{ex}\label{ex:1984}
In the Introduction we mentioned that computing the 
subspace of Shimura equivalent half-integral weight
forms is necessary
for applying Waldspurger's Theorem~\cite{Waldspurger} to a given
integral weight form. In this example we shall
illustrate this, working with the weight $2$ newform $F$
corresponding to the elliptic curve
\[
E/\Q \; :\quad Y^2 = X^3 + X + 1.
 \]
The elliptic curve $E$ has conductor $N_0=496 = 16\times 31$ and 
$j$-invariant $6912/31$; in particular $E$ does not have 
complex multiplication.
Let $F\in S_2^{\mathrm{new}}(496)$ be the corresponding newform
with trivial character given by the Modularity Theorem; 
$F$ has the following $q$-expansion,
\[ 
F= q - 3q^5 + 3q^7 - 3q^9 - 2q^{11} - 4q^{13} - q^{19} + O(q^{20}).
 \]

In order to apply Waldspurger's Theorem~\cite{Waldspurger}
we would like to find a suitable $N$ divisible by
$2 N_0=992$, and character $\chi$, such that the
summand
$S_{3/2}(N,\chi,F)$ is non-trivial, and we also
need to compute an eigenbasis for this summand. 
Let $\chi_0$ be the trivial 
character. Using Theorem~\ref{thm:algo} we found that 
$S_{3/2}(992,\chi_0,F)= \{0\}$. Next we considered level $N=1984$. 
The space $S_{3/2}(1984,\chi_0)$ 
is $119$-dimensional (see \cite{C-O}). 
From Theorem~\ref{thm:basisS0} it follows that the space $S_0(1984,\chi_0)$ 
is $1$-dimensional and is spanned by the theta series given by 
$\sum_{m=1}^{\infty} \chi_{-1}(m)mq^{m^2}$
where $\chi_{-1}$ is the Dirichlet character of conductor $4$ given by $\chi_{-1}(\cdot)= \kro{-1}{\cdot}$. 
Using our algorithm we obtain the following decomposition of $S_{3/2}^\prime(1984,\chi_0)$
where the $F_i$ vary over the 
newforms of levels dividing $992$ and trivial character, deg $\Q(F_i)$ denotes the degree of number field generated by the coefficients of $F_i$ 
and $D_i$ is the dimension of the space $S_{3/2}(1984,\chi_0,F_i)$.
In this table we group conjugate newforms on the same row.
\begin{table}[htbp]
\begin{center}
\begin{tabular}{cc}
\begin{minipage}{0.5\hsize}
\begin{center}
\begin{tabular}{|c|c|c|c|}
\hline
Level & Newforms $F_i$ & deg $\Q(F_i)$ & $D_i$\\ \hline
31 & $F_1$, $F_2$ & 2 & 12 \\
32 & $F_3$ & 1 & 0 \\
62 & $F_4$ & 1 & 9 \\
62 & $F_5$, $F_6$ & 2 & 9 \\
124 & $F_7$& 1 & 6 \\
124 & $F_8$& 1 & 6 \\
248 & $F_9$& 1 & 3 \\
248 & $F_{10}$ & 1 & 3 \\
248 & $F_{11}$ & 1 & 3 \\
248 & $F_{12}$, $F_{13}$ & 2 & 3 \\
248 & $F_{14},\dots,F_{16}$ & 3 & 3 \\
496 & $F_{17}$ & 1 & 3 \\
496 & $F_{18}$ & 1 & 3 \\
496 & $F_{19}$ & 1 & 3 \\
496 & $F_{20}$ & 1 & 2 \\
\hline
\end{tabular}
\end{center}
\end{minipage}
\begin{minipage}{0.5\hsize}
\begin{center}
\begin{tabular}{|c|c|c|c|}
\hline
Level & Newforms $F_i$ & deg $\Q(F_i)$ & $D_i$\\ \hline
496 & $F_{21}$ & 1 & 2 \\
496 & $F_{22}$ & 1 & 1 \\
496 & $F_{23}$, $F_{24}$ & 2 & 3 \\
496 & $F_{25}$, $F_{26}$ & 2 & 0 \\
496 & $F_{27}$, $F_{28}$ & 2 & 1 \\
496 & $F_{29},\dots,F_{31}$ & 3 & 3 \\
992 & $F_{32}$, $F_{33}$ & 2 & 0 \\
992 & $F_{34}$, $F_{35}$ & 2 & 0 \\
992 & $F_{36},\dots,F_{38}$ & 3 & 0 \\
992 & $F_{39},\dots,F_{41}$ & 3 & 0 \\
992 & $F_{42},\dots,F_{45}$ & 4 & 0 \\
992 & $F_{46},\dots,F_{49}$ & 4 & 0 \\
992 & $F_{50},\dots,F_{55}$& 6 & 0 \\
992 & $F_{56},\dots,F_{61}$& 6 & 0 \\
    & &  &    \\
\hline
\end{tabular}
\end{center}
\end{minipage}
\end{tabular}
\end{center}
\end{table} 

Our $F$ is in fact $F_{17}$ in the above table.
From the table we see that the space 
$S_{3/2}(1984,\chi_0,F_{17})$ is $3$-dimensional.
Our algorithm also gives an eigenbasis  
$\{f_1,f_2,f_3\}$ given by the following $q$-expansions:
\begin{equation*}
\begin{split}
f_1 & = q^3 + q^{43} - 2q^{75} + 2q^{83} + q^{91} + 3q^{115} - 3q^{123} + O(q^{145}):= \sum_{n=1}^{\infty} a_nq^n\\
f_2 & = q^{15} + q^{23} - q^{31} + 2q^{55} + q^{79} - 3q^{119} + O(q^{145}):= \sum_{n=1}^{\infty} b_nq^n\\
f_3 & = q^{17} + q^{57} + q^{65} + 2q^{73} - q^{89} - q^{105} + q^{137} + O(q^{145}):= \sum_{n=1}^{\infty} c_nq^n.
\end{split}
\end{equation*}

Using Waldspurger now we can prove the following statement. One can find the details of the proof in \cite{SomaIII}.

\begin{prop}
Let $f=f_1+f_2+\sqrt{2} f_3=\sum d_n q^n$. For positive square-free $n \equiv 1$, $3$, $7 \pmod{8}$,
\begin{enumerate}
 \item [(a)] \[
\mathrm{L}(E_{-n},1) = \frac{2^{(\nu_{31}(n)+1)} \Omega_{E_{-1}}}{\sqrt{n}}
\cdot d_n^2.\] 
 \item[(b)] Let $E_{-n}$ has rank zero. Then assuming the Birch and Swinnerton-Dyer conjecture, 
\[
|\SHA{E_{-n}/\Q}|= \frac{2^{(\nu_{31}(n)+1)}}{\prod_{p}{c_p}} \cdot d_n^2
\]
where the Tamagawa numbers $c_p$ of $E_{-n}$ are given by
\[
c_2=\begin{cases} 
1 & n \equiv 3,7 \pmod{8}\\
2 & n \equiv 1,5 \pmod{8},
\end{cases}
\qquad 
c_{31}=\begin{cases}
1 & 31 \nmid n, \\
4 & 31 \mid n, \kro{n/31}{31}=1\\
2 & 31 \mid n, \kro{n/31}{31}=-1,
\end{cases}
\]
and $c_p=\# E_{-1}(\F_p)[2]$ for $p \mid n$, $p \ne 31$, and $c_p=1$ for all other primes $p$.
\item[(c)] Suppose $\kro{n}{31} = -1$. Then assuming the Birch and Swinnerton-Dyer Conjecture,
\[\mathrm{Rank}(E_{-n}) \geq 2 \Leftrightarrow d_n = 0.\]
\end{enumerate}
\end{prop}

\end{ex}

\end{document}